\documentclass{ajour} 
\usepackage{amssymb,array,delarray,latexsym} 

\newlength{\sh}
\newlength{\jmr}
\newlength{\jfc}
\newlength{\bernd}
\newlength{\jones}
\newlength{\mati}
\newlength{\sil}
\newlength{\koi}
\newtheorem{kho}{Khovanski's Theorem on Real Fewnomials}

\newtheorem{cor}{Corollary}

\newtheorem{main}{Main Theorem} 
\newtheorem{smale}{Smale's Theorem} 

\newtheorem{plaisted}{Plaisted's Theorem} 

\newtheorem{thm}[main]{Theorem}
\newtheorem{prob}{Problem}
\newtheorem{rem}{Remark}	
\newtheorem{ex}{Example}

\newcommand{\pp}{\bbf{P}} 
\newcommand{\np}{\bbf{NP}} 
\newcommand{\bpp}{\bbf{BPP}} 
\newcommand{\bp}{\bbf{BP}}

\newcommand{\eps}{\varepsilon}
\newcommand{\cA}{\mathcal{A}}
\newcommand{\cD}{\mathcal{D}}

\newcommand{\cO}{\mathcal{O}}

\newcommand{\conv}{\mathrm{Conv}}

\newcommand{\thth}{{\underline{\mathrm{th}}}} 

\newcommand{\st}{ {\underline{ \mathrm{st} } }  }

\newcommand{\Q}{\mathbb{Q}}
\newcommand{\R}{\mathbb{R}}
\newcommand{\C}{\mathbb{C}}
\newcommand{\N}{\mathbb{N}}
\newcommand{\Z}{\mathbb{Z}}

\newcommand{\Zn}{\Z^n}

\newcommand{\Rn}{\R^n}

\newcommand{\Rsn}{{(\R^*)}^n}
\newcommand{\Csn}{{(\C^*)}^n}
\renewcommand{\qed}{$\blacksquare$}

\newcommand{\cC}{\mathcal{C}}

\newcommand{\bO}{\bbf{O}}
\newcommand{\vol}{\mathrm{Vol}}

\begin{document}

\title{Some Speed-Ups and Speed Limits for Real Algebraic Geometry } 
\titlerunninghead{Speed-Ups and Speed Limits for Real Algebraic Geometry } 

\dedication{In memory of Gian-Carlo Rota.} 

\author{J.\ Maurice Rojas\thanks{ This research was partially funded by Hong 
Kong UGC Grant \#9040402. } }  
\authorrunninghead{J.\ Maurice Rojas}  

\affil{Department of Mathematics, 
City University of Hong Kong, 
83 Tat Chee Avenue, 
Kowloon, HONG KONG. }
\email{ mamrojas@math.cityu.edu.hk,  
Web-page: http://www.cityu.edu.hk/ma/staff/rojas}  


\date{\today} 

\abstract{ 
We give new positive and negative results, some conditional, on 
speeding up computational algebraic geometry over the reals: 
\begin{enumerate} 
\item{ A new and sharper upper bound on the number of connected 
components of a semi-algebraic set. Our bound is novel in that it is 
stated in terms of the volumes of certain polytopes and, for a large class 
of inputs, beats the best previous bounds by a factor 
exponential in the number of variables. } 
\item{ A new algorithm for approximating the real roots of certain sparse 
polynomial systems. Two features of our algorithm are (a) arithmetic 
complexity {\bf polylogarithmic} in the degree of the underlying 
complex variety (as opposed to the super-linear dependence in earlier 
algorithms) and (b) a simple and efficient generalization to certain 
univariate exponential sums. } 
\item{ Detecting whether a real algebraic surface (given as the 
common zero set of some input straight-line programs) is not smooth can 
be done in polynomial time within the classical Turing model (resp.\ BSS 
model over $\C$) only if $\pp\!=\!\np$ (resp.\ $\np\!\subseteq\!\bpp$). } 
\end{enumerate} 
The last result follows easily from an unpublished observation of Steve 
Smale. } 

\section{Introduction and Main Results}
\label{sec:intro} 
We provide new speed-ups for some fundamental computations in real algebraic 
geometry. Our techniques are motivated by recent results from algebraic 
geometry but the proofs are almost completely elementary. We then conclude 
with a discussion of how much farther these techniques can still be pushed. 

In particular, we significantly improve the best previous upper 
bounds on the number of connected components of a {\bf 
semi-algebraic\footnote{A {\bf semi-algebraic set} is simply a 
subset of $\Rn$ defined by the solutions of a finite collection of 
polynomial inequalities.} set}, and we 
exhibit a new class of polynomial systems over the 
real numbers which can be solved within polylogarithmic time. As 
for complexity lower bounds, we show that if singularity  
detection for curves over $\C$ can be done in polynomial time then, 
depending on the computational model, we must 
have $\pp\!=\!\np$ or $\np\!\subseteq\!\bpp$. This can also be thought of as a 
lower bound on the complexity of elimination theory, and 
immediately implies an analogous result on singularity detection 
for real algebraic surfaces.  

This work is a part of an ongoing program by the author 
\cite{rio,real,gcp,dio} to dramatically sharpen current complexity bounds 
from algebraic geometry in terms of more intrinsic geometric invariants. 
We will give precise statements of these results shortly, so 
let us begin by considering the number of connected components of a 
semi-algebraic set.

\subsection{Sharper Intrinsic Bounds}\mbox{}\\  
\label{sec:bound} 
The topology of semi-algebraic sets is intimately related 
to complexity theory in many ways. For example, the seminal 
work of Dobkin, Lipton, Steele, and Yao \cite{dl79,sy82} (see 
also \cite[Ch.\ 16]{bcss}) relates upper bounds on the number of connected 
components to lower bounds on the algebraic circuit complexity of 
certain problems. More directly, upper bounds on connected components 
are an important ingredient in complexity upper bounds for 
the first order theory of the reals \cite{bpr}. 

Our first main theorem significantly improves earlier bounds on 
the number of connected components by Oleinik, 
Petrovsky, Milnor, Thom, and Basu \cite{op,milnor,thom,basu}.\footnote{These 
papers actually bound the sum of the {\bf Betti numbers}, which in turn is an  
upper bound on the number of connected components. Our bounds 
can be extended to bound the sum of the Betti numbers as well, but this 
extension will be addressed in future work. } 
The main novelty of our new bound is its greater sensitivity to the 
monomial term structure of the input polynomials. Letting $\bO$ and 
$\hat{e}_i$ respectively denote the origin and the $i^{\thth}$ standard 
basis vector in $\R^N$, $x\!:=\!(x_1,\ldots,x_n)$, and normalizing 
$k$-dimensional volume $\vol_k(\cdot)$ 
so that the standard $k$-simplex $\Delta_k\!:=\!\{x\!\in\!\R^k \; | \; 
x_1,\ldots,x_n\!\geq\!0 , \sum_j x_j\!\leq\!1\}$ has volume $1$, our result is 
the following. 
\begin{main}\mbox{}\\  
\label{main:betti}
\label{MAIN:BETTI}
\noindent 
Let $f_1,\ldots,f_{p+s}\!\in\!\R[x_1,\ldots,x_n]$ and suppose 
$S\!\subseteq\Rn$ is the solution set of the following 
collection of polynomial inequalities: 
\begin{eqnarray*} 
f_i(x)\!&\!=\!&\!0, \ \ \ i\!\in\!\{1,\ldots,p\} \\ 
f_{p+i}(x)\!&\!>\!&\!0, \ \ \ i\!\in\!\{1,\ldots,s\} 
\end{eqnarray*} 
Let $Q\!\subset\!\Rn$ be the convex hull of the 
union of $\{\bO,\hat{e}_1,\ldots,\hat{e}_n\}$ 
and the set of all $a$ with $x^a\!:=\!x^{a_1}_1\cdots x^{a_n}_n$ a  
monomial term of some $f_i$. Then $S$ has at most 
\[ \mbox{}\hspace{-.3cm}\min\{n+1,\frac{s+1}{s-1}\}2^ns^n\vol_n(Q) 
(\mathrm{for \ } s\!>\!0) \ \ \mathrm{or}   \ \ 2^{n-1}\vol_n(Q) 
(\mathrm{for \ } s\!=\!0) \]  
connected components. 
\end{main} 
\noindent 

In section \ref{sec:ex} we show that this bound is at least 
as good as (and frequently much better than) the aforementioned earlier 
bounds.  Our bound also considerably simplifies, and is competitive with, 
an earlier polytopal bound of Benedetti, Loeser, and Risler 
\cite[Prop.\ 3.6]{blr}.
(We note that their polytopal bound, in addition to some minor 
restrictions on the $f_i$, applies only when $s\!=\!0$ and $p\!\leq\!n$.)  

It is interesting to note that there are sharper (even optimal) upper bounds 
relating polytope volumes and connected components for 
{\bf complex} varieties, beginning with the remarkable 
work of Bernshtein, Kushnirenko, and Khovanski \cite{bkk} 
a bit over twenty years ago. (See also \cite{dk}.\footnote{ 
We also point out that the classical B\'ezout's theorem \cite{mumford} is 
optimal only for a small class of polynomial systems. So the results of 
\cite{bkk} include B\'ezout's theorem as a very special case.}) 
However, as far as the author is aware, Main Theorem \ref{main:betti} presents 
the first nontrivial general upper bounds on the number of connected 
components of {\bf semi-algebraic} sets with this combinatorial flavor.   
The work of Benedetti, Loeser, and Risler \cite{blr} appears to be the 
first occurence of polytopal bounds for the case where $s\!=\!0$ and 
$p\!\leq\!n$ (i.e., 
certain real algebraic sets). 

\begin{rem} 
Finding an {\bf optimal} upper bound on the number of connected 
components of a semi-algebraic set, even in the 
special case of nondegenerate real algebraic sets, remains an open problem. 
\qed 
\end{rem} 

Our bound can be further improved in various ways and 
this is detailed in section \ref{sec:betti}. In particular, 
we give sharper versions tailored for certain special cases (e.g., 
compact hypersurfaces and real algebraic sets), and 
we prove analogues (for {\bf all} our bounds) depending only on $n$, $s$, and 
the number of monomial terms which appear in at least one $f_i$. 
Khovanski appears to have been the first to consider bounds of 
this type for the case where $s\!=\!0$ and $p\!\leq\!n$ \cite{few}. 

The techniques involved in our proof of Main Theorem \ref{main:betti}, when 
combined with other recent results of the author \cite{gcp}, 
also yield similar improvements on the complexity of 
quantifier elimination over real-closed fields. 
This will be pursued in a forthcoming paper of the author. 

\subsection{Superfast Real Solving for Certain Fewnomial Systems}\mbox{}\\ 
The complexity of solving systems of {\bf fewnomials} 
(polynomials with few monomial terms\footnote{Results on fewnomials usually 
hold on a much broader class of functions: the so-called {\bf Pfaffian} 
functions \cite{few}.}) has been addressed only recently. Indeed, 
the vast majority of work in 
computational algebra has so far been stated only in terms of degrees 
of polynomials, thus ignoring the finer monomial term structure. Notable 
exceptions include \cite{cks} (solving a single univariate fewnomial over 
$\Z$ in polynomial time), \cite{lenstra} (solving a single univariate 
fewnomial over $\Q$ in polynomial time), and \cite{real,mp98,gcp,gls99} 
(solving polynomial systems over $\R$ or $\C$ within time near polynomial 
in the degree of the underlying complex variety).  

While it is more or less intuitively clear what it means to solve 
a polynomial system over $\Z$ or $\Q$, let us state a motivating 
problem to clarify what we mean by solving over $\R$: 
\begin{prob}
\label{prob:1}
Can one $\eps$-approximate all the roots of a univariate fewnomial of degree 
$d$, within the interval $[0,R]$, using significantly less than 
$\Theta(d\log\log\frac{R}{\eps})$ arithmetic steps? \qed 
\end{prob} 

In particular, an important alternative statement is the following: 
\begin{prob} 
\label{prob:2} 
Can the complexity of solving fewnomials be {\bf sub-linear} in the degree of 
the underlying complex variety? \qed 
\end{prob} 
Finding such ``super-fast'' algorithms is nontrivial, even for {\bf bi}nomials 
(i.e., quickly finding $d^{\thth}$ roots) \cite{ye}. The asymptotic complexity 
limit stated in Problem \ref{prob:1}, up to a factor polylogarithmic in $d$, 
is the best current bound for solving a 
general univariate polynomial of degree $d$ over $\C$ \cite{neffreif}. In 
particular, the existence of faster algorithms for finding just the {\bf 
real} roots of a degree $d$ fewnomial was unknown until now. 

Our next main theorem gives an affirmative answer to Problem \ref{prob:2}, for 
certain fewnomial systems and univariate exponential sums over $\R$. More 
precisely, if $f(x)\!=\!\sum_{a\in \cA} c_ax^a$, where $\cA\!\subset\!\R$ 
is finite and the coefficients $c_a$ are all real, we call $f$ a {\bf 
(real) exponential $\bbf{k}$-sum}. When $\cA\!\subset\!\Z$, we define the 
{\bf degree} of such an $f$ to be $\max_{a,a'\in \cA}\{a-a'\}$. Otherwise, 
we set $\bbf{\deg(f)}\!:=\!\max\{a-a'\}/\min\{1,\min\{a-a'\}\}$, where the 
second minimum ranges over all {\bf distinct}\footnote{ 
We declare the degree of any monomial to be $0$.}  
$a,a'\!\in\!\cA$. We also say that $f$ has {\bf $\bbf{j}$ sign 
alternations} iff there are $j$ distinct pairs $(a,a')\!\in\!\cA^2$ such that 
$c_ac_{a'}\!<\!0$, $\cA\cap (a,a')\!=\!\emptyset$, and $a'\!>\!a$. So,  
for instance, $47x^{2.53}-10.3x^{0.9}-\pi-10x^{-3}-x^{-5.5}$ has just one 
sign alternation but $x^3-2x+2$ has two. Finally, when 
$\cA\!\subset\!\Z$, we simply call $f$ a {\bf $\bbf{k}$-nomial}. 
\begin{main} 
\label{main:tri} 
\label{MAIN:TRI} 
Let $f$ be any exponential $k$-sum of degree $d$ with at most one sign 
alternation. Then, given an oracle for evaluating $x^r$ for any $x,r\!\in\!\R$, 
one can $\eps$-approximate all the roots\footnote{...and of course count their 
number} of $f$ in $(0,R)$ 
using $\cO(k(\log d+\log\log\frac{R}{\eps}))$ arithmetic operations 
over $\R$ (including oracle calls). In particular, restricting to $k$-nomials 
and removing the oracle, we can still do the same using 
$\cO(k\log d(\log d+\log\log\frac{R}{\eps}))$ arithmetic 
operations over $\R$, with $d$ agreeing with the usual  
degree of a univariate Laurent polynomial. 
\end{main} 

We point out that even the {\bf tri}nomial case is difficult. For 
example, while one can count the number of real roots of a trinomial of the 
form $x^d+ax+b$ within $\cO(\log d)$ arithmetic operations \cite{richard} 
(regardless of sign alternations), doing the same for general trinomials 
was an open problem until recently \cite{ry}. 
Also, even from a numerical point of 
view, the use of Newton's method is subtle for trinomials: 
It is known that deciding whether a given initial point converges  
to a root of $x^3-2x+2$ is undecidable in the BSS model over $\R$ 
(see \cite[Sec.\ 2.4]{bcss} and \cite{barna}). 
Nevertheless, this need not stop us from finding {\bf some} good 
starting point, as we will soon see. 

Our algorithm, aside from an algebraic trick, closely follows an 
algorithm of Ye \cite{ye} (for a particular class of analytic 
functions) which efficiently blends binary search  
and Newton's method. By combining these ideas with a few facts 
on the {\bf Smith normal form} of an integral matrix \cite{unimod}, 
we can also derive the following complexity result on 
binomial systems. 
\begin{main} 
\label{main:bi} 
\label{MAIN:BI} 
Let $c_1,\ldots,c_n\!\in\!\R\!\setminus\!\{0\}$ and let $[d_{ij}]$
be any $n\times n$ matrix with nonnegative integer entries. Finally, let
$f_i:=x^{d_{i1}}_1\cdots x^{d_{in}}_n+c_i$ for all $i$. Then we can
$\eps$-approximate all the roots of $f_1\!=\cdots=\!f_n\!=0$ in the {\bf
orthant wedge} $\{x\!\in\!\Rn \; | \; x_1,\ldots ,x_n\!\geq\!0 , 
\sum_i x^2_i\!\leq\!R^2\}$ within
\[ \cO((n+\log \max|d_{ij}|)^{6.376}) \mathrm{ \ bit \ operations,} \]
followed by
\[\cO\left(\log \left|\det[d_{ij}]\right| \left[n^3\log^2(n \max|d_{ij}|)+
\log\log\frac{R}{\eps}\right]\right) \] 
rational operations over $\R$. 
\end{main} 
If the above binomial system has only finitely many {\bf complex} roots,
then their number is exactly $|\det[d_{ij}]|$. This follows easily from
{\bf Bernshtein's theorem} \cite{bkk}. It is also interesting
to note that the fastest previous general
(sequential) algorithms for polynomial system solving over $\R$ or $\C$, 
when applied to binomial systems, run in time {\bf polynomial} in 
$|\det[d_{ij}]|$ \cite{mp98,gcp,gls99} --- that is, super-linear in the degree 
of the underlying complex variety. 

One can of course solve slightly more general systems of fewnomials 
by threading together the algorithms of Main Theorems \ref{main:tri} and 
\ref{main:bi}. We will say more on the likelihood of farther-reaching 
extensions of our last two results after first discussing a 
result relating complexity classes and singularities. 

\begin{rem} 
Finding $\eps$-approximations of roots within a suitable region is far 
from the strongest notion of solving a polynomial system. In particular,   
the spacing between roots, which of course dictates the $\eps$ one 
should choose, must be taken into account. A more complete and 
elegant framework would be to include the {\bf condition number} 
\cite{bcss} of the input fewnomial system in all complexity bounds. 
It is thus the author's intent that the preceding fewnomial complexity bounds 
be interpreted as a first step in this direction. \qed 
\end{rem} 

\subsection{Obstructions to Superfast Degeneracy Detection}\mbox{}\\ 
The preceding two algorithmic results circumvent degeneracy 
problems in simple but subtle ways. For instance, Main Theorem 
\ref{main:tri} clearly deals with equations having at most 
one positive real root, while the binomial systems of Main Theorem 
\ref{main:bi} are easily seen to have no repeated complex roots (cf.\ 
section \ref{sec:few}). Thus, the respective hypotheses of these results 
(restricting sign alternations and/or number of monomial terms) allow 
us to approximate roots without stopping for a singularity check. 

It seems hard to completely solve a system of equations without 
knowing something about its degeneracies, either a priori or 
during run-time. So let us present a result which 
gives solid evidence that detecting degeneracies may be quite 
difficult. In what follows, unless otherwise mentioned, we use 
the standard {\bf sparse encoding} 
for multivariate polynomials \cite{plaisted,hnam}. Thus the {\bf size} of a 
polynomial like $x^{d}+x-47$ will be $\Theta(\log d)$ and not $\Theta(d)$, 
whether in the Turing model or the BSS model over $\C$. 
\begin{main} 
\label{main:steve} 
\label{MAIN:STEVE} 
Suppose any of the following problems can be solved in 
polynomial time via a Turing machine (resp.\ BSS machine over 
$\C$). Then $\pp\!=\!\np$ (resp.\ $\np\!\subseteq\!\bpp$). 
\begin{enumerate}
\item{Decide if an input polynomial $f\!\in\!\Z[x_1]$ (resp.\ 
$f\!\in\!\C[x_1]$) vanishes at a $d^{\thth}$ root of unity, 
where $d\!=\!\deg(f)$. } 
\item{Decide if two input polynomials $f,g\!\in\!\Z[x_1]$ (resp.\ 
$f,g\!\in\!\C[x_1]$) have a common root. }
\item{Given a nonzero input polynomial $f\!\in\!\Z[x_1,x_2]$ (resp.\  
$f\!\in\!\C[x_1,x_2]$) decide if the curve $\{(x_1,x_2)\!\in\!(\C^*)^2 \; | 
\; f(x_1,x_2)\!=\!0\}$ has a singularity. } 
\item{Given input polynomials $f,g\!\in\!\Z[x_1,x_2,x_3,x_4]$ (resp.\  
$f,g\!\in\!\R[x_1,x_2,x_3,x_4]$), in the straight-line program encoding, 
defining a surface $S\!\subset\!\R^4$, decide 
if $S$ has a singularity. } 
\item{Given any finite subset $\cA\!\subset\!\Z^2$ and a 
vector of coefficients $(c_a \; | \; a\!\in\!\cA)\!\in\!\Z^{\#\cA}$ 
(resp.\ $\in\!\!\C^{\#\cA}$), decide if the $\cA$-discriminant 
of the bivariate polynomial $\sum_{a\in\cA} c_ax^a$ vanishes.}  
\end{enumerate} 
\end{main} 
\begin{rem} 
Note that in problem (4) we are already given that $S$ is a surface. 
Determining whether this is true or not turns out to be $\np$-hard 
(resp.\ $\np_\R$-complete) in the Turing model (resp.\ BSS 
model over $\R$) \cite{npdim}. 
\qed 
\end{rem} 
\noindent 
For any $\cA\!\subset\!\Zn$, the {\bf $\bbf{\cA}$-discriminant}, $\cD_A$, 
is defined to be the unique (up to sign) irreducible polynomial in $\Z[c_a 
\; | \; a\!\in\!\cA]$ such that $f_\cA(x)\!:=\!\sum_{a\in \cA} c_ax^a$ has a 
singularity in its zero set (in $(\C^*)^n) \Longrightarrow \cD_A\!=\!0$ 
\cite{gkz94}. This important operator lies at the heart of {\bf sparse 
elimination theory}, which is the part of algebraic geometry surrounding this 
paper. 

The $\cA$-discriminant in fact contains all known 
multivariate resultants and discriminants as special cases, and 
also appears in residue theory and hypergeometric functions \cite{gkz94}. 
Thus, a corollary of our last main result is that sparse elimination 
theory, even in low dimensions, might lie beyond the reach of $\pp$. 
\begin{rem}
It is interesting to note that nontrivial lower bounds on the complexity of 
computing $\cA$-discriminants in the {\bf one}-dimensional case 
$\cA\!\subset\!\Z$ are unknown. However, it is easy to show (via 
\cite[pg.\ 274]{gkz94}) that one can at 
least find $D_\cA$ in {\bf polynomial time} when $\cA\!\subset\!\Zn$ has less 
than $n+3$ elements. \qed  
\end{rem}  

We will prove our main theorems in order of appearance, but 
first let us return to our study of semi-algebraic sets  
to see some examples. 

\section{Comparing Upper Bounds on the Number of Connected Components}
\label{sec:ex}
Here we briefly compare our first main theorem to earlier 
bounds on the number of connected components of a 
semi-algebraic set. 

In summary, we can compare our new bound to earlier bounds (stated 
in terms of total degree) in very simple polyhedral terms: Let $\Delta_Q$ 
denote the smallest 
scaled standard $n$-simplex, $\gamma\Delta_n$, containing $Q$. Then, since 
volume is monotonic under containment, our bounds are least favorable when 
$Q\!=\!\Delta_Q$. However, in practice it will frequently be the 
case that $Q$ has much smaller volume that $\Delta_Q$, thus accounting 
for improvements as good as a factor exponential in $n$.  

\subsection{At Least One Inequality} 
Assume $s\!>\!0$ temporarily. Letting $d$ be the maximum of the total 
degrees of the $f_i$, the best previous general upper 
bounds, quoted from 
\cite[Ch.\ 16, Prop.\ 5]{bcss} and \cite{basu} respectively, were 
$(sd+1)(2sd+1)^n$ and $(p+s)^n\cO(d)^n$. (The first bound 
is an improved version of a bound due to Milnor, Oleinik, Petovsky, 
and Thom \cite{op,milnor,thom}.) Our bound is no worse 
than $\min\{n+1,\frac{s+1}{s-1}\}(2sd)^n$ (better than both 
preceding bounds) and is frequently much better. 
Consider the following examples: 
\begin{ex} {\bf (Spikes)} 
Suppose we pick all the $f_i$ to have the same monomial term 
structure, and in such a way that $Q$ has small volume but 
great length some chosen direction. In particular, let us 
assume that the only monomial terms occuring in the $f_i$ are 
$1,x_1,\ldots,x_{n-1}$ and $(x_1\cdots x_n)$, $(x_1\cdots x_n)^2$, $\ldots$, 
$(x_1\cdots x_n)^D$. Then 
it is easy to check that $Q$ is a ``long and skinny'' bypyramid, with one 
apex at the origin and the other at $(D,\ldots,D)\!\in\!\Rn$. 
We then obtain, via two simple determinants, that 
$\vol_n(Q)\!=\!D+1$ and thus our bound reduces to 
$\min\{n+1,\frac{s+1}{s-1}\}2^ns^n(D+1)$. 
However, the aforementioned older bounds are easily seen to reduce to 
$(nsD+1)(2nsD+1)^n$ and $((p+s)\cO(nD))^n$. \qed  
\end{ex} 
\begin{ex} {\bf (Bounded Multidegree)} 
Suppose now that instead of bounding the total degree of the 
$f_i$, we only require that the degree of $f_i$ with respect 
to any $x_j$ be at most $d'$. It is then easy to check that 
$Q$ is an axes parallel hypercube with side length $d'$. 
So our new bound reduces to $\min\{n+1,\frac{s+1}{s-1}\}(2sd')^n$. However, 
the old bounds are easily seen to reduce to $(snd'+1)(2snd'+1)^n$ 
and $((p+s)\cO(nd'))^n$. \qed 
\end{ex} 

\subsection{Real Algebraic Sets} 
Assume now that $s\!=\!0$. Then the aforementioned earlier upper 
bounds respectively reduce to $d(2d-1)^n$ and $(p\cO(d))^n$. Specializing 
Main Theorem \ref{main:betti}, we obtain a bound which is no worse than 
$2^{n-1}d^n$ (neglibly worse than the first, better than the second), and is 
frequently much better. This can easily be 
seen by reconsidering our last two examples in the case $s\!=\!0$. 
(We leave this as an exercise.)
 
However, let us now make a fairer comparison to another polytopal bound --- 
that of Benedetti, Loeser, and Risler \cite[Prop.\ 3.6]{blr}. 
\begin{rem} 
The bound \cite[Prop.\ 3.6]{blr} was published with several 
typographical errors. Following inquiries from the author, 
Francois Loeser kindly responded via three e-mails with the  
following corrections:\footnote{Professor Loeser states that 
these corrections were also checked with Jean-Jacques Risler, 
one of the other authors of \cite{blr}.} in the notation of their bound, 
a hypothesis of $k\!\leq\!n$ was missing. Also, 
in part (a) of their statement, the quantity $\Phi(\Delta)$  
should be replaced by $\theta^n_k(\Delta)$, and the last 
sum should be replaced by the main quantity from Prop.\ 3.1. 
Finally, in part (c), all j's should be capitalized, and 
$\theta$ should be replaced by $\theta^n_{k-\#J}$. \qed 
\end{rem} 

The bound \cite[Prop.\ 3.6]{blr} 
has a recursive definition based on {\bf mixed} volumes 
\cite{gk94,mvcomplex}. For the sake of brevity, we will focus on the four  
examples given in \cite{blr}. 

\begin{ex} {\bf (Four Examples from \cite{blr})} 
Examples (A), (B), (C), and (D) of \cite[Sec.\ 4]{blr} concern 
polynomial systems of the following shape: 
(A) $c_0+c_1x^a+c_2y^b$ (one polynomial, two variables), 
(B) $c_0+c_1x^{a_1}_1+\cdots+c_nx^{a_n}_n$ (one polynomial, 
$n$ variables), (C) $c_0+c_1x+c_2y+c_3(xy)^a$ (one polynomial, 
two variables), and (D) $(c_0+c_1x^a+c_2y^b,c_3+c_4x^b+c_5y^b+c_6(xy)^b)$ 
(two polynomials, two variables), where the $c_i$ are 
real constants and $a,b\!\in\!\N$. 

The polytopal bound of \cite{blr}, when applied to these examples 
in the above order, respectively evaluates to $2ab+4$, $2a_1\cdots a_n+
\mathrm{Lower \ Order \ Terms}$, $8a$, and $2ab-b^2+\mathrm{Lower \ Order \ 
Terms}$. None of the preceding lower order terms is stated 
explicitly in \cite{blr}, and it appears that the last value is incorrect. 
However, a closer examination of their (corrected) bound 
respectively yields $2ab+4$, $2(a_1+2)\cdots(a_n+2)$, $8a$, and 
$8b^2+6ab+8$. 

Main Theorem \ref{main:betti} is easily seen to 
respectively evaluate to $2ab$, $2^{n-1}a_1\cdots a_n$, $4a$, and $4ab$ 
for these examples. \qed 
\end{ex} 

More generally, it is not hard to check that 
our bound is usually better than that of \cite{blr} when 
$n$ is small or $p$ is close to $n$. (Indeed, the bound 
of \cite{blr} does not cover the case $p\!>\!n$.)  
However, the bound from \cite{blr} usually wins 
when $p$ is a small constant and $n$ is large. The author 
hopes to combine the techniques here with those of 
\cite{blr} in future work. 
 
\section{\mbox{Proving Main Theorem \ref{main:betti}}} 
\label{sec:betti} 
We will first prove a sharper version of Main Theorem \ref{main:betti} for 
compact hypersurfaces, and then successively generalize to 
the case of real algebraic and semi-algebraic sets. 
Along the way, we give analogues of our upper bounds depending only 
on $n$, $s$, and the number of monomial terms. 

\begin{rem} 
Throughout this section, ``nonsingular'' (or ``smooth'') for a real algebraic 
variety will mean that the underlying {\bf complex} variety is nonsingular in 
the sense of the usual Jacobian criterion (see, e.g., \cite{mumford}). \qed 
\end{rem} 

\subsection{Point-Free Compact Zero Sets of a Single Polynomial} 
We begin with the following important special case of 
Main Theorem \ref{main:betti}. This lemma is also frequently 
significantly sharper than many earlier results and may be of 
independent interest.  
\begin{lemma}
\label{lemma:compact} 
Following the notation of Main Theorem \ref{main:betti}, 
suppose $p\!=\!1$, $s\!=\!0$, and $S$ is 
compact but has no zero-dimensional components. Then 
$S$ has at most $\frac{1}{\min\{2,n\}}\vol_n(Q')$ connected components, where 
$Q'$ is the convex hull of the union of $\{\bO\}$ and the set of all $a$ 
with $x^a$ a monomial term of $f_1$. 
\end{lemma}
\noindent 
{\bf Proof:} The main idea will be to show that (for $n\!\geq\!2$) the number 
of connected components is bounded above by half the number of 
critical points of a projection of a perturbed version of 
$S$. This idea is quite old, but we will introduce an unusual 
projection which permits a much sharper upper bound than 
before. The case $n\!=\!1$ of our bound is trivial, so let us 
assume $n\!\geq\!2$ henceforth. 

Consider $\tilde{f}\!:=\!f_1+\delta$, for some $\delta\!\in\!\R$ 
to be selected later. By Sard's theorem \cite{hirsch}, there is a 
set $W\!\subseteq\!\R$ of full measure such that $\delta\!\in\!W 
\Longrightarrow S_\delta\!=\!\{ 
x\!\in\!\Rn \; | \; \tilde{f}\!=\!0\}$ is nonsingular (and a 
hypersurface). 
Also, via a simple homotopy argument, $S$ and $S_\delta$ are both compact and 
have the same number of connected 
components, for $|\delta|$ sufficiently small. (Much stronger versions of this 
fact can be found in \cite{basu}.) Furthermore, note that for all but 
finitely many $\delta$, no connected component of $S_\delta$ lies  
inside the union of the coordinate hyperplanes. We will pick $\delta\!\neq\!0$ 
so that all these conditions, and one more to be described below, hold. 
  
Now consider the function $x^a$, with 
$a\!\in\!\Zn\!\setminus\!\{\bO\}$ 
to be selected later. Clearly, any connected component of 
$S$ (not lying in a hypersurface of the form $x^a\!=\!\mathrm{constant}$) 
must have at least two special points: one locally maximizing, and the other 
locally minimizing, $x^a$. Since there are only finitely many connected 
components (by any earlier bound, e.g., \cite{op}), and every 
component contains a curve, there must therefore be 
an $a\!\in\!\Zn\!\setminus\!\{\bO\}$ so that every component 
(not lying entirely within the union of coordinate hyperplanes) 
contributes at least two critical points of $x^a$. Pick $a$ 
in this way, subject to the additional minor restricition 
that the g.c.d.\ of the coordinates of $a$ is $1$. 

Note that the critical points of the function $x^a$ on $S_\delta$ are 
just the solutions in $\Rn$ of 
\[ (\star) \ \ \ \ \ \tilde{f}\!=\!\frac{\partial 
\tilde{f}}{\partial y_2}\!=\cdots =\!\frac{ \partial \tilde{f}}{\partial 
y_n}\!=\!0, \]  
where the $y_i$ are new variables to be described shortly.  
Our final condition on $\delta$ (which is easily seen to hold for all but 
finitely many $\delta$) will simply be that all real solutions to the 
above polynomial system lie in $\Rsn\!:=\!(\R\!\setminus\!\{0\})^n$. 
Note also that a corollary of all our assumptions so far is that the 
number of {\bf complex} solutions of ($\star$) is finite. (This 
follows immediately from Sard's theorem, and the fact that 
the complex solutions of ($\star$) form an algebraic set.) 

We are now essentially done: The number of connected components of 
$S$ and $S_\delta$ are the same, and the latter quantity is bounded 
above by half the number of critical points (on $S_\delta$) of the function 
$x^a$. This number of critical points can be computed in terms 
of polytope volumes as follows: Via the Smith normal form \cite{smith}, 
we can find an invertible change of variables on $\Rsn$ such that 
$y_1\!:=\!x^a$ and $y_2,\ldots,y_n$ are monomials in the $x_i$. 
Furthermore, this change of variables induces the action of a unimodular 
matrix on the exponent vectors of $\tilde{f}$. In particular, 
$\tilde{f}$ can be considered as a polynomial in 
$\R[y^{\pm 1}_1,\ldots,y^{\pm 1}_n]$ and the 
number of monomial terms (and Newton polytope volume) of $\tilde{f}$ 
is preserved under this change of variables. Thus, up to a monomial change of 
variables, the critical points of the function $x^a$ on $S_\delta$ are exactly 
the solutions in $\Rsn$ of ($\star$). 

The key to our new bound is to finish things off by picking a bound {\bf 
other} than B\'ezout's theorem here. In particular, by Bernshtein's 
theorem \cite{bkk}, 
the number of solutions in $\Csn$ is at most the {\bf mixed volume} of $Q'$ 
and $n-1$ other polytopes with translates contained in $Q'$. 
By the monotonicity of the mixed volume \cite{buza}, 
the latter quantity is at most the mixed volume of 
$n$ copies of $Q'$ and, by the definition of mixed volume, 
this is just $\vol_n(Q')$. \qed 

We point out that a key ingredient in our proof is 
that the monomial change of variables we use (as opposed to the  
linear changes of variables used in most earlier 
treatments) preserves sparsity. This allows us to take full advantage of more 
powerful and refined techniques to bound the number of real roots, 
and thus get new bounds on the number of real connected components. 
For example, substituting Bernshtein's theorem for B\'ezout's theorem  
in the older proofs would not have yielded any significant improvement. 

However, we need not have been so heavy-handed and only used tools over $\C$. 
We could have also used the following alternative bound on the number of real 
roots. 
\begin{kho}\mbox{}\\
{\bf (Special Case)} \cite[Sec.\ 3.12, Cor.\ 6]{few}
Suppose that for all $i\!\in\!\{1,\ldots,n\}$, 
$f_i\!\in\!\R[x_1,\ldots,x_n,m_1,\ldots,m_k]$ has total degree $q_i$, where 
the $m_j$ are monomials in $x$.  Let $k$ denote the number of monomial terms 
which appear in at least one of $f_1,\ldots,f_n$. 
Assume further that the variety $S$ defined 
by $f_1,\ldots,f_n$ is zero-dimensional and nonsingular. Then $S$ has 
at most $(1+\sum_i q_i)^k 2^{k(k-1)/2} \prod q_i$ connected components in the 
positive orthant. \qed 
\end{kho}
We call any set of the form $\{x\!\in\!\Rn \; | \; \pm x_1,\ldots,\pm x_n
\!\geq\!0\}$ a {\bf closed orthant}. When all signs are positive we call 
the corresponding closed orthant the {\bf nonnegative orthant}. The analogous 
constructions where all inequalities are strict are, respectively, an {\bf 
open} orthant and the {\bf positive} orthant. 

As an immediate corollary, our proof above yields the following 
alternative upper bound on the number of components of a 
smooth compact real algebraic hypersurface. 
\begin{cor} 
\label{cor:kho} 
Following the notation of lemma \ref{lemma:compact}, assume 
further that $S$ is a smooth compact hypersurface. Then 
the number of connected components of $S$ is at 
most $2^{n-1}(n+1)^{k+1}2^{k(k+1)/2}$. In particular, 
$S$ has at most $\frac{1}{2}(n+1)^k2^{k(k-1)/2}$ connected components 
contained entirely within the positive orthant. 
\end{cor} 
\noindent 
{\bf Proof:} 
Following the notation of our last proof, 
note that multiplying any equation of ($\star$) by a monomial 
in $y_1,\ldots,y_n$ does not affect the roots in $\Rsn$. Thus, we can 
assume ($\star$) has only $k+1$ distinct monomial terms. 
Also note that the monomial change of variables $x\mapsto y$  
maps orthants onto orthants, and that the case $n\!=\!1$ 
is trivial. The first portion of our corollary then 
follows immediately from our last proof (using Khovanski's Theorem 
on Fewnomials with $q_1\!=\cdots=\!q_n\!=1$ instead of Bernshtein's Theorem), 
upon counting roots in all open orthants. The second portion follows even more 
easily, upon observing that we do not need $\delta$ 
if we only want to count critical points in an open orthant. \qed 

\subsection{The Case of Real Algebraic Varieties} 
The next step in proving Main Theorem \ref{main:betti} is to 
increase the number of polynomials allowed and drop the 
compactness hypothesis. Again, the following result is 
frequently much sharper than many earlier bounds and may also be of 
independent interest. 
\begin{lemma}
\label{lemma:variety} 
Following the notation of Main Theorem \ref{main:betti}, 
suppose now that $s\!=\!0$, so that $S$ is a real algebraic 
variety, not necessarily smooth or compact. Then 
$S$ has at most $2^{n-1}\vol_n(Q)$ connected components.  
\end{lemma}
\noindent 
{\bf Proof:} The main trick is to reduce to the 
case considered by our preceding lemma. In particular, 
define $F_{\delta,\eps}\!:=\!f^2_1+\cdots+f^2_p+\eps^2(\sum_i x^2_i)-
\delta^2\!\in\!\R[x_1,\ldots,x_n]$ and let $S_{\delta,\eps}$ be the set of 
real zeroes of $F_{\delta,\eps}$. 
It then follows that for sufficiently small (and suitably restricted) 
$\delta,\eps\!>\!0$, $S_{\delta,\eps}$ is a smooth compact 
hypersurface and the number of connected components of $S_{\delta,\eps}$ is 
no smaller than the number of connected components of $S$. The proof 
of this fact is standard and a very clear account can be found 
in \cite[Sec.\ 16.1]{bcss}. 

In any event, the number of connected components of $S_{\delta,\eps}$ 
is clearly at most $\frac{1}{2}\vol_n(\conv(2Q'\cup \{2\hat{e}_1,
\ldots,2\hat{e}_n\}))$, thanks to our preceding lemma. Since 
the last quantity is just $\frac{1}{2}\cdot2^n\vol_n(Q)$ we are done. \qed 

We can combine the proof of lemma \ref{lemma:variety} with Khovanski's 
Theorem on Fewnomials to obtain the following generalization of corollary 
\ref{cor:kho}. This result, while giving a slightly looser bound than an 
earlier result of Khovanski \cite[Sec.\ 3.14, Cor.\ 5]{few}, removes all 
the nondegeneracy assumptions from his result. 
\begin{cor} 
\label{cor:me} 
Following the notation and assumptions of lemma \ref{lemma:variety}, 
the number of connected components of $S$ is also bounded 
above by\\ $4^{n-\frac{1}{2}}(2n+1)^{k+1}2^{k(k+1)/2}$.
\end{cor} 

\noindent 
{\bf Proof:} Combining the proofs of lemmata \ref{lemma:variety} 
and \ref{lemma:compact}, and since we are only counting roots in 
$\Rsn$, we see that the number of 
connected components is at most half the number of 
solutions in $\Rsn$ of the following polynomial system: 
\[ (\star\star) \ \ \ \bar{F}_{\delta,\eps}\!=\!y_2\frac{\partial 
\bar{F}_{\delta,\eps}} {\partial y_2}\!=\cdots =\!y_n\frac{\partial 
\bar{F}_{\delta,\eps}} {\partial y_n}\!=\!0,\] 
where $\bar{F}_{\delta,\eps}$ is the variant of $F_{\delta,\eps}$ 
where we substitute $\sum_i y^2_i$ for $\sum_i x^2_i$. (It is 
a simple exercise to verify that the proof of lemma \ref{lemma:variety} 
still goes through with this variation.) Now simply note, via the chain rule 
of calculus, that every polynomial in ($\star\star$) is of degree at most $2$ 
in $y_1,\ldots,y_n$ and the set of monomials appearing in 
$f_1,\ldots,f_p$. Also note that the polynomials in ($\star\star$) 
are polynomials in a total of $k+1$ monomial terms. So by Khovanski's Theorem 
on Real Fewnomials, and counting roots in all open orthants, we are done. \qed 

\subsection{Extending to Semi-Algebraic Sets} 
We are now ready to prove Main Theorem \ref{main:betti}. 

\noindent
{\bf Proof of Main Theorem \ref{main:betti}:} 
We reduce again, this time to lemma \ref{lemma:variety}. 
The trick here is to note that every connected component 
of $S$ is in turn a connected component of $S'$ where 
$S'\!:=\!\{x\!\in\!\Rn \; | \; f_1(x)\!=\cdots =\!f_p(x)\!=\!0, 
f_{p+1}(x)\!\neq\!0, \ldots,f_{p+s}(x)\!\neq\!0 \}$. Every 
connected component of $S'$ is in turn a projection 
(onto the first $n$ coordinates) of a connected component 
of $S''$, where $S''\!\subset\!\R^{n+1}$ is the real zero set of 
the polynomial system $(f_1,\ldots,f_p,-1+z\prod^{p+s}_{i=p+1} 
f_i)$. This reduction is not new and appears, among other places, 
in \cite[Sec.\ 16.3]{bcss}. 

Now lemma \ref{lemma:variety} tells us that 
the number of connected components of $S''$ is at most 
$2^n$ times the $(n+1)$-dimensional volume of $\conv(P_1\cup 
(P_2\times\hat{e}_{n+1}))$, 
where $P_1$ (resp.\ $P_2$) is the union of $\{\bO,\hat{e}_1,\ldots,
\hat{e}_n\}$ and the Newton polytopes 
of $f_1,\ldots,f_p$ (resp.\ the {\bf Minkowski sum} of the 
Newton polytopes of $f_{p+1},\ldots,f_{p+s}$). However, it is 
a simple exercise to show that $P_2\!\subseteq\!P_3$ where $P_3$ is the 
union of $\{\bO,\hat{e}_1,\ldots,\hat{e}_n\}$ and the Newton polytopes of 
$f_{p+1},\ldots,f_{p+s}$, scaled by a factor of $s$. Now note that 
$P_2\!\subseteq\!Q$, $P_3\!\subseteq\!sQ$ 
and $\conv(P_1\cup (P_2\times\hat{e}_{n+1}))\!\subseteq\!
\conv(Q\cup (sQ\times\hat{e}_{n+1}))$. 

If $s\!>\!1$ then the last 
polytope is in turn contained in a pyramid $P$ with apex 
at $(0,\ldots,0,\frac{-1}{s-1})$ and base $Q\times\hat{e}_{n+1}$. 
So we obtain that the number of connected components 
of $S$ is at most $2^n\vol_{n+1}(P)\!=\!2^n\frac{s+1}{s-1}
\vol_n(sQ)\!\!=\!\!\frac{s+1}{s-1}2^ns^n\vol_n(Q)$. 

If $s\!=\!1$ then $\conv(Q\cup (sQ\times\hat{e}_{n+1}))\!=\![\bO,\hat{e}_{n+1}]
\times Q$. So, similar to the previous case, the number of connected 
components of $S$ is at most $2^n\vol_{n+1}(P)\!=\!2^n n\vol_n(Q)$. 

Now note that the number of connected components of $S$ will always be 
at most $\min\{n+1,\frac{s+1}{s-1}\}2^ns^n\vol_n(Q)$, with the possible 
exception of the case $(n,s)\!=\!(1,2)$. So we need 
only check this final case. However, this is almost 
trivial, separating the cases $p\!>\!0$ and $p\!=\!0$. \qed 

We can give an alternative version of Main Theorem \ref{main:betti}, solely 
in terms of $n$, $s$, and $k$, as follows. 
\begin{thm} 
Following the notation and assumptions of Main Theorem \ref{main:betti}, 
the number of connected components of $S$ is also bounded 
above by $4^{n-\frac{1}{2}}(s+1)^n(2(n+1)(s+1)+1)^{k+1}2^{k(k+1)/2}$. \qed 
\end{thm} 

The proof is very similar to that of corollary \ref{cor:me}, save 
only that we substitute the polynomial system from the 
proof of Main Theorem \ref{main:betti} into the construction of 
$\bar{F}_{\delta,\eps}$. In particular, we eventually obtain a 
system of $n+1$ polynomials of degree $2(s+1)$ in a total of $k+1$ 
monomials, thus allowing yet another application of Khovanksi's 
beautiful theorem on fewnomials.  

\section{Alpha Theory and Proving Main Theorems \ref{main:tri} and 
\ref{main:bi}} 
\label{sec:few}
The proof of Main Theorem \ref{main:tri} hinges on 
{\bf gamma theory} \cite{bcss}, which gives useful 
criteria for when Newton's method converges quadratically. In particular, 
we will need the following elementary analytic lemma. 
\begin{lemma} 
For any monotonic function $\phi : \R \longrightarrow \R$, let 
$\gamma_\phi$ satisfy $\sup_{k>1}|\frac{\phi^{(k)}(x)}{k!\phi'(x)}|^{
\frac{1}{k-1}}\!\leq\!\frac{\gamma_\phi}{x}$. Then, 
for $\phi(x)\!=\!x^r$, we may take $\gamma_\phi$ equal to 
$\lceil |r| \rceil$, $2$ or $1$, according as $r\!\in\!(-\infty,-1)\cup 
(1,\infty)$, $r\!\in\!(0,1)$, or $r\!\in\!(-1,0)$. More 
generally, if $\phi\!=\!\phi_1+\phi_2$ with $\phi_1$ and 
$\phi_2$ both convex and either both increasing or both 
decreasing, then we can take $\gamma_\phi\!=\!\max\{\gamma_{\phi_1},
\gamma_{\phi_2}\}$. \qed 
\end{lemma}
\noindent  
The first part is a simple exercise while the second part is a proposition 
from \cite{ye}. 

We are now ready to sketch the proof of Main Theorem \ref{main:tri}.

\noindent
{\bf Proof of Main Theorem \ref{main:tri}:} 
We begin by changing our function $f$ slightly. First let 
$M$ be largest exponent occuring in the $k$-sum $f$ and let 
$m$ be the smallest real number so that $x^m$ is a monomial 
term of $f$ with {\bf positive} coefficient. (We assume, 
by multiplying by $-1$ if necessary, that the leading coefficient 
of $f$ is positive.) By dividing out by $x^m$ we may assume that $m\!=\!0$. 
Via the change of variables $x\!=\!y^{1/M}$, we may further assume that 
$M\!=\!1$. In particular, we now obtain that $f$ is a 
sum of two increasing convex functions: one a positive 
linear combination of powers of $x$ (with exponents 
in $(0,1]$), the other a negative 
linear combination of powers of $x$ (with exponents in $(-\infty,0)$). 

By our preceding lemma, we may take $\gamma_f\!=\!d$ (the degree 
of $f$) since $d$ is no smaller than the degree of our original $f$. We now 
invoke the hybrid algorithm from \cite[Theorem 3]{ye}: 
This algorithm allows us to $\eps$-approximate the real roots 
of $f$ in $(0,R)$ using $\cO(\log \gamma_f + \log\log\frac{R}{\eps})\!=\!
\cO(\log d + \log\log\frac{R}{\eps})$ function evaluations and 
arithmetic operations. To conclude the first part of this main theorem, 
inverting the change of variables we made requires another 
$\cO(\log d + \log\log\frac{R}{\eps})$ operations via the same algorithm 
(since taking $n^{\thth}$ roots is the same as solving an exponential 
$2$-sum). However, 
we may have decreased the accuracy of our $\eps$-approximation. 
So we just begin by solving to accuracy $\min\{\eps^{M-m},\eps\}$ instead to 
obtain the first part of our main theorem. (Note also that 
evaluating $f$ requires $k$ uses of our oracle.) 

To obtain the second part of our theorem, we simply use the 
same algorithm without the oracle. This simply introduces another 
factor of $\log d$ since monomials can now be evaluated by the usual 
repeated squaring trick. \qed 
 
Main Theorem \ref{main:bi} only needs a special case of 
Main Theorem \ref{main:tri}. In fact, \cite{ye} contains a slightly 
modified algorithm for the binomial case with an even better complexity 
bound of $\cO(\log d \log \log \frac{R}{\eps})$, which we will use 
below. However, we will also require some refined quantitative 
facts about the Smith normal form of a matrix. 
\begin{lemma} \cite{unimod} 
\label{lemma:unimod} 
Let $A\!=\![a_{ij}]$ be any $n\times n$ matrix with entries only in $\Z$ and 
define $h_A$ to be $\log(2n+\max|a_{ij}|)$. Then, within 
$\cO^*((n+h_A)^{6.375})$ bit operations, one can find matrices 
$U$, $D$, $V$ with the following properties:
\begin{enumerate} 
\item{$U$ and $V$ both have determinant $\pm 1$ and entries only in $\Z$.} 
\item{$D$ is diagonal and has entries only in $\Z$.} 
\item{$UAV\!=\!D$} 
\item{$\det A$ is the product of the diagonal elements of $D$ and 
$h_U,h_V\!=\!\cO(n^3(h_A+\log n)^2)$. } 
\end{enumerate} 
\end{lemma} 

\noindent 
{\bf Proof of Main Theorem \ref{main:bi}:} 
We begin by immediately applying the Smith normal form to 
our matrix $[d_{ij}]$. (This accounts for the bit 
operation count.) Clearly then, we have reduced 
to the case of $n$ binomials of the form 
$x^{d_1}_1-\gamma_1,\dots,x^{d_n}_n-\gamma_n$.  
The real roots of this  polynomial system can then be 
$\eps$-approximated by $n$ applications of Main Theorem \ref{main:tri}. 
Since $\sum_i \log d_i\!=\!\log \prod_i d_i\!=|\det [d_{ij}]|$, this 
accounts for almost all of the second bound. 

To conclude, note that we must still invert our change of 
variables. By lemma \ref{lemma:unimod}, computing this 
monomial map is almost the final contribution to our second 
complexity bound. The only missing part is the fact that we 
may have needed  more accuracy at the beginning of our 
algorithm. Lemma \ref{lemma:unimod} also tells us how much more accuracy we 
need, thus finally accounting for all of our second complexity bound. \qed 

\section{Smale's Theorem and Main Theorem \ref{main:steve}} 

\label{sec:smale} 
We begin with the following result of Plaisted. 
\begin{plaisted} 
\cite{plaisted} 
Deciding if an input polynomial $f\!\in\!\Z[x_1]$  
coefficients) vanishes at a $d^{\thth}$ root of unity, where 
$d\!=\!\deg(f)$, is $\np$-hard. 
\qed 
\end{plaisted} 
\noindent 
In the above (and in what follows) $f$ is given in the sparse encoding, 
so coefficients {\bf and} exponents are measured by bit-length. 

The following unpublished result of Steve Smale gives an 
intriguing extension of Plaisted's result via computations 
over new rings. 
\begin{smale}
Suppose we can decide, within polynomial time  
relative to the BSS model over $\C$, if an input polynomial 
$f\!\in\!\C[x_1]$ vanishes at a $d^{\thth}$ root of unity, 
where $d\!=\!\deg(f)$. Then $\np\!\subseteq\!\bpp$. 
\end{smale} 

\noindent
{\bf Proof:} Given any complexity class $\cC$ over the Turing model, 
consider its extension $\cC_\C$ to the BSS model over $\C$. 
It is then a simple fact that $\cC$ is contained in the {\bf Boolean part} 
of $\cC_\C$, $\bp(\cC_\C)$ \cite{bpp}. However, we will make use of 
an inclusion going the opposite way: $\bp(\cC_\C)\!\subseteq\!\cC^{\bpp}$ 
\cite{bpp}. Applying this to the problem at hand, we thus see that the 
hypothesis of Smale's theorem, thanks to Plaisted's Theorem, implies that 
$\np\!\subseteq\!\bp(\pp_\C)\!=\!\pp^{\bpp}\!=\!\bpp$. So we are done. \qed 

Our final main theorem then follows from some simple 
reductions to problem (1) from the statement. 

\noindent 
{\bf Proof of Main Theorem \ref{main:steve}:}  
First note that the assertion concerning 
problem (1) follows immediately from Smale's Theorem and 
Plaisted's Theorem. It thus suffices to successively reduce (1) to 
special cases of all the other problems. 

The assertion for (2) is then clear, since via the 
special case $g(x)\!=\!x^d-1$, any polynomial time 
algorithm for (2) would give a polynomial time algorithm 
for (1). 

On the other hand, a polynomial time algorithm for problem (5) 
would imply a polynomial time algorithm for problem (2). 
This is because problem (2) is essentially the decision problem 
of whether the {\bf sparse resultant} of $f$ and $g$ \cite{gkz94} is zero. 
Via the {\bf Cayley trick} \cite{gkz94}, the $\cA$-discriminant 
for $\cA\!=\!P\cup(Q\times \hat{e}_2)$ (where $P$ and 
$Q$ are respectively the {\bf supports}\footnote{The support of a 
polynomial is simply the set of its exponent vectors (fixing 
an ordering on the variables).}  of $f$ and $g$) is 
exactly the sparse resultant of $f$ and $g$. So this portion is done. 

Note also that (3) is just a reformulation of (5). 

As for (4), via the Jacobian criterion for singularities \cite{mumford} 
applied to the real and imaginary parts of the input to (3), a 
polynomial time algorithm for (4) (using the straight-line program encoding 
for the input) would immediately imply a polynomial time algorithm for (3) 
(using the straight-line program encoding for the input). Such an 
algorithm would then immediately be a polynomial time 
algorithm for (3) with inputs given in the sparse encoding. \qed 

\section{Acknowledgements} 
The author thanks Felipe Cucker, Askold Khovanski, and Steve Smale for some 
very useful discussions. In particular, he thanks Felipe Cucker for 
pointing out an elegant proof of Smale's Theorem. Special thanks 
also go to Francois Loeser for pointing out the excellent paper 
\cite{blr}. 

\footnotesize
\bibliographystyle{acm}

\end{document}